\documentclass{amsart}
\newtheorem{theorem}{Theorem}[section]
\newtheorem{lemma}[theorem]{Lemma}

\newtheorem{definition}[theorem]{Definition}

\numberwithin{equation}{section}

\def\C{\mathbb C}
\def\R{\mathbb R}
\def\X{\mathbb X}

\def\Y{\mathbb Y}

\def\cal{\mathcal}
\def\cD{\cal D}
\def\tD{\tilde{{\cal D}}}
\def\F{\cal F}
\def\tf{\tilde{f}}
\def\tg{\tilde{g}}

\begin{document}
\title[Katznelson-Tzafriri Theorems]{Katznelson-Tzafriri Type Theorems for Individual Solutions of Evolution Equations}
\author{Nguyen Van Minh}
\address{Department of Mathematics, University of West Georgia, Carrollton, GA 30118}
\email{vnguyen@westga.edu}
\thanks{The author thanks the referee for reading carefully the
manuscript and for making useful remarks.}

\date{\today}
\begin{abstract}In this paper we present an extension of the
Katznelson-Tzafriri Theorem to the asymptotic behavior of
individual solutions of evolution equations $u'(t) =Au(t)+f(t)$.
The obtained results do not require the uniform continuity of
solutions as well as the well-posedness of the equations. The
method of study is based on a recently developed approach to the
spectral theory of functions that is direct and free of
$C_0$-semigroups.
\end{abstract}
\keywords{Katznelson-Tzafriri Type Theorem, reduced spectrum of a
function, asymptotic behavior} \subjclass{34G10; 47D06}

\maketitle


\section{Introduction and Statement of Results}
This is a companion paper of \cite{min} in which the author
developed a new approach to the spectral theory of functions on
the line and on the half line to study the asymptotic behavior of
solutions of evolution equations of the form
\begin{equation}\label{eq}
\dot u(t) = Au(t) +f(t),
\end{equation}
where $t$ is in $J$ ($J$ is either $\R^+$ or $\R$), $A$ is a
closed linear operator on a Banach space $\X$, $f\in BC(J,\X)$.
This spectral theory is simple, free of $C_0$-semigroup theory, so
it can apply to large classes of solutions and equations with
general conditions.

\bigskip
In this paper we will discuss further applications of the theory
by proving continuous analogs of the two well known Katznelson -
Tzafriri Theorems in \cite{kattza} to individual bounded solutions
of evolution equations of the form (\ref{eq}) on the half line.
Our main results are stated in the two theorems below (Theorems
\ref{the 1}, \ref{the 2}). For this purpose we introduce classes
of function spaces $\F$ as subspaces of $BC(\R^+,\X)$ that satisfy
so-called {\it Condition $F^+$ or $F^{++}$} (see the definition in
the next section). As an example of such function spaces $\F$ one
may take $C_0(\R^+,\X)$.

\begin{theorem}\label{the 1}
Let $\F$ be a function space that satisfies Condition $\F^+$, $A$
be a linear operator on a Banach space $\X$ such that $\sigma (A)
\cap i\R \subset \{ 0\}$, and let $u\in BC(\R^+,\X)$ be a
classical solution of (\ref{eq}) on $\R^+$ with $f\in \F$ such
that $Au(\cdot )\in BC(\R^+,\X)$. Then,
\begin{equation}
Au(\cdot ) \in \F .
\end{equation}
\end{theorem}

Recall that a function $g\in L^1(\R) $ is of spectral synthesis
with respect to a closed subset $E \subset \R$ if there exists a
sequence $(g_n)$ in $L^1(\R)$ such that
\begin{enumerate}
\item $lim_{n\to\infty} \|  g - g_n\|_1 = 0$; \item each of the
Fourier transforms $\F g_n$ vanishes on a neighborhood of $E$.
\end{enumerate}
\begin{theorem}\label{the 2}
Let $\F$ be a function space that satisfies Condition $\F^{++}$,
 $u\in BC(\R^+,\X)$ be a mild solution of (\ref{eq}) on
$\R^+$ with $f\in\F$, and let $\phi\in L^1(\R^+)$ be of spectral
synthesis with respect to $i\sigma (A)\cap \R$. Then, $w\in \F$,
where
\begin{equation}
w(t):=  \int^\infty _0 \phi (s) u(t+s)ds .
\end{equation}
\end{theorem}
The proofs of these theorems will be given in Section 3. Notations
and a short introduction into previous results in \cite{min} will
be given in the next section. For the works related to Theorem
\ref{the 1} we refer the reader to \cite{arebathieneu,arepru}. Our
Theorem \ref{the 1} is an individual version of the part "(ii)
$\Rightarrow$ (i)" of \cite[Theorem 3.10]{arepru}. In fact, if $A$
generates an eventually differentiable bounded $C_0$-semigroups
$T(t)$, $\F$ is chosen to be $C_0(\R^+,\X)$ and $f=0$, then since
$T(\tau )\X \subset D(A)$ for some $\tau \ge 0$, for each $x\in
\X$ the function $u(t):= T(t)x$ is a classical solution (with
$t\ge \tau$). Therefore, our Theorem \ref{the 1} applies. Notice
that under the assumptions of Theorem \ref{the 1}, a version of
"(i) $\Rightarrow$ (ii)" of \cite[Theorem 3.10]{arepru} does not
make sense in our context because there may not be bounded
solutions to an ill-posed equation. As for Theorem \ref{the 2},
our result extends the individual version of Katznelson-Tzafriri
Theorem given in \cite[Theorem 5.1]{batneerab} that applies to
homogeneous equations generating bounded $C_0$-semigroups. In
fact, if $A$ generates a bounded $C_0$-semigroup $T(t)$, $\F$ is
chosen to be $C_0(\R^+,\X)$ and $f=0$, then for each $x\in \X$,
and $\phi \in L^1(\R^+)$
\begin{eqnarray*}
T(t) \int^\infty _{-\infty} \phi (\xi )T(\xi )x d\xi &=&
\int^\infty
_{-\infty} \phi (\xi )T(t+\xi )x d\xi \\
&=& \int^\infty _{-\infty} \phi (\xi )u(t+\xi ) d\xi := w(t),
\end{eqnarray*}
where $u(t):= T(t)x$. Our result shows in particular that
$\lim_{t\to\infty} w(t)=0.$

\bigskip
For more complete accounts of related works we refer the reader to
the monographs \cite{arebathieneu,nee} and their references, as well
as the papers
\cite{arepru,baspry,batneerab,chitom,eststrzou,mus,vu}.

\section{Preliminaries}
Throughout the paper we denote by $\R$ the real line, by $\R^+$
the positive half line $[0,\infty )$, by $R^-$ the negative half
line $(-\infty ,0]$, and by $\X$ a Banach space over the complex
plane $\C$. If $A$ is a linear operator on a Banach space $\X$,
$D(A)$ stands for its domain; $\sigma(A)$ stands for its spectrum.
$L(\X)$ stands for the Banach space of all bounded linear
operators in $\X$ with the usual norm $\| \cdot \|$. If
$\lambda\in \rho (A)$, then $R(\lambda ,A)$ denotes the resolvent
$(\lambda -A)^{-1}$. In this paper we also use the following
notations:
\begin{enumerate}
\item $BC(J,\X)$ is the space of all $\X$-valued bounded and
continuous functions on $J$, where $J$ is either $\R$, or $\R^+$;
\item $L^\infty (J,\X)$ and $L^1(J,\X)$ are the space  of all
$\X$-valued measurable and essentially bounded functions on $J$,
and the space  of all $\X$-valued measurable and integrable
functions on $J$, respectively; \item $C_0(\R^+,\X):= \{ f\in
BC(\R^+,\X): \ \lim_{t\to \infty} f(t)=0\}$; \item If $A$ is a
linear operator on $\X$, then the operator of multiplication by
$A$ on $BC(J,\X)$, denoted by ${\cal A}$, is defined on $D({\cal
A}):= \{ g\in BC(J,\X):\ g(t)\in D(A), \ \mbox{for all} \ t\in J,
Ag(\cdot )\in BC(J,\X)\}$, by ${\cal A}g=Ag(\cdot )$ for each
$g\in D({\cal A})$.
\end{enumerate}

\begin{definition}
An $\X$-valued continuous function  $u$ on $\R^+$ is said to be a
mild solution of (\ref{eq}) on $\R^+$ (with $f\in L^\infty (\R^+
,\X)$) if for every $t\in \R^+$, $\int^t_0 u(s)ds \in D(A)$, and
\begin{equation}\label{mild sol}
u(t) -u(0) = A\int^t_0 u(s)ds +\int^t_0 f(s)ds , \quad \ \mbox{for
all} \ t\in \R^+ .
\end{equation}
A continuously differentiable function $u$ on $\R^+$ is said to be
a classical solution of (\ref{eq}) on $\R^+$ (with $f\in BC (\R^+
,\X)$) if for every $t\in \R^+$, $u(t) \in D(A)$, and
\begin{equation}\label{class sol}
u'(t)  = Au(t) +f(t) , \quad \ \mbox{for all} \ t\in \R^+ .
\end{equation}
\end{definition}

\medskip
Let us introduce the so-called  {\it Conditions $F^+$ and
$F^{++}$} for a function space ${\cal F} \subset BC(\R^+,\X)$. Let
$\psi$ be any function in $L^\infty(\R^+, \X)$ (or $L^\infty(\R^-,
\X)$, respectively) throughout the paper we will identify it with
its natural extension to $\psi \in L^\infty(\R, \X)$ by setting
$\psi (t)=0$ for all $t<0$ (or $t>0$, respectively).
\begin{definition}
 A function space ${\cal F} \subset BC(\R^+,\X)$ is
said to satisfy {\it Condition $F^{+}$} if
\begin{enumerate}
\item It is closed, and contains $C_0(\R^+,\X)$; \item If $g\in
{\cal F}$, then the function $\R^+\ni t\mapsto e^{i\xi t}g(t)\in
\X$ is in ${\cal F}$ for all $\xi\in\R$; \item For each $h\in
{\cal F}$, $Re\lambda >0, \ Re\eta <0$, the function $y(\cdot),
z(\cdot )$, defined as
\begin{equation}
y(t) = \int^\infty _t e^{\lambda (t-s)} h(s)ds, \ \ z(t) =
\int^t_0 e^{\eta (t-s)} h(s)ds ,\quad  t\in \R^+
\end{equation}
are in ${\cal F}$;  \item For each $B\in L(\X)$ and $f\in {\cal F}$,
the function $Bf(\cdot )$ is in $\cal F$.
\end{enumerate}
If in addition to the above conditions, ${\mathcal F}$ satisfies the
following:
\begin{enumerate}
\item[v)] For each function
$\psi \in L^1(\R^- )$, $\psi * g\in {\mathcal F}$ for each $g\in
{\mathcal F}$,
\end{enumerate}
then ${\mathcal F}$ is said to satisfy {\it Condition} $F^{++}$.
\end{definition}
As an example of a function space that satisfies Conditions $F^+$
and $F^{++}$, we can take ${\cal F}=C_0(\R^+,\X)$. Another function
space that satisfies Condition $F^{++}$ is $AA(\R^+,\X)$, the space
of all restrictions to $\R^+$ of the $\X$-valued almost automorphic
functions. Note that $AA(\R^+,\X)$ contains non-uniformly continuous
functions, so it is not a subspace of $BUC(\R^+,\X)$.

\bigskip
Consider the quotient space $\Y := BC(\R^+,\X)/ {\cal F}$ whose
elements are denoted by $\tilde{g}$ with $g\in BC(\R^+,\X)$. We will
use $\tilde{\cD}$ to denote the operator induced by the differential
operator $\cD$ on $\Y$ which is defined as follows: The domain
$D(\tD )$ is the set of all classes that contains a differentiable
function $g\in BC(\R^+,\X)$ such that $g'\in BC(\R^+,\X)$; $\tD \tg
:= \tilde{g '}$ for each $\tg \in D(\tD )$.

\bigskip
It is easy to see that for $f\in BC(\R^+,\X)$ (see \cite[page
13]{min})
\begin{eqnarray}\label{resolvent of tilde D}
R(\lambda ,\tD )\tf  (t) &=& \begin{cases} \int^\infty _t
e^{\lambda
(t-s)} \tf (s)ds  ,\quad  Re\lambda > 0, \ t\in \R^+ ,\\ ~ \\
- \int^t_0 e^{\lambda (t-s)} \tf (s)ds  \quad Re\lambda <0, \ \
t\in \R^+.
\end{cases}\\
&=&\begin{cases} \int^\infty_0 e^{-\lambda \eta} \tf (t+\eta )d\eta ,\quad  Re\lambda > 0, \ t\in \R^+ ,\\ ~ \\
-\int^t_0 e^{\lambda \eta}\tf (t-\eta ) d\eta , \quad Re\lambda
<0, \ \ t\in \R^+.\nonumber
\end{cases}
\end{eqnarray}

Under the above notations, the operator $\tilde{\cD}$ is a closed
operator with $\sigma (\tilde{\cD} ) \subset i\R$.

\begin{definition}
Let ${\cal F}$ be a function space that satisfies Condition
$F^{++}$, and let $f\in BC(\R^+,\X)$. Then the {\it reduced
spectrum} of $f$ with respect to ${\cal F}$, denoted by $sp^+_{{\cal
F}}(f)$, is defined to be the set of all reals $\xi \in\R$ such that
$R(\lambda , \tilde{\cD})\tilde{f}$, as a complex function of
$\lambda$ in $\C \backslash i\R$, has no holomorphic extension to
any neighborhood of $i\xi$ in the complex plane.
\end{definition}
Since $sp^+_{{\cal F}}(f)$ is the same for all elements $f$ in a
class $\tg$, the use of the notation $sp^+_{{\cal F}}(\tg)$ makes
sense. The following theorem was proved in \cite{min}:

\begin{theorem}\label{the 3}
Let ${\cal F}$ be a function space of $BC(\R^+,\X)$ that satisfies
Condition $F^{+}$, and let $f$ be in $BC(\R^+,\X)$ such that
$sp^+_{\cal F}(f)$ is countable. Moreover, assume that
\begin{equation}
\lim_{\alpha \downarrow 0} \alpha R(\alpha +i\xi ,\tD )\tf =0
\end{equation}
for all $\xi \in sp^+_{\cal F}(f)$. Then, $f\in {\cal F}$.
\end{theorem}

For each function $f\in L^\infty (\R,\X)$, $\xi\in\R$,
$Re\lambda\not= 0$, let us recall that
\begin{eqnarray}\label{2.6}
R(\lambda,{\cal D})f (\xi )&=& \begin{cases}
\begin{array}{ll}
\int^  {\infty}_0  e^{-\lambda \eta }f(\xi +\eta )d\eta
&(\mbox{if}\ Re\lambda > 0)\\ \\
-\int _ {-\infty}^0  e^{-\lambda \eta}f(\xi +\eta )d\eta     &
(\mbox{if } \ Re\lambda < 0).
\end{array} \end{cases}
\end{eqnarray}
Recall also that the Carleman transform of $f\in L^\infty (\R,\X)$
is defined as
\begin{eqnarray}
\hat f (\lambda )&=& \begin{cases}
\begin{array}{ll}
\int^  {\infty}_0  e^{-\lambda \eta }f( \eta )d\eta
&(\mbox{if}\ Re\lambda > 0)\\ \\
-\int _ {-\infty}^0  e^{-\lambda \eta}f( \eta )d\eta     &
(\mbox{if } \ Re\lambda < 0).
\end{array} \end{cases}
\end{eqnarray}
It can be shown in the same manner as in the proof of
\cite[Proposition 2.3]{liunguminvu} that the set of $\zeta \in \R$
such that $\hat f (\lambda )$ has no holomorphic extension to any
neighborhood of $i\zeta$ (Carleman spectrum of $f$) coincides with
the set of all $\zeta\in\R$ such that $R(\lambda,{\cal D})f$ has
no holomorphic extension to any neighborhood of $i\zeta$. We will
denote this set by $sp(f)$ and call it the spectrum of $f$. The
reader is referred to \cite{arebathieneu,liunguminvu} for more
information on this concept of spectrum.

\section{Proof of the Main Results}
The main idea of proving Theorems \ref{the 1}, and \ref{the 2} is to
apply Theorem \ref{the 3}. We will fix $\F$ as a function space that
satisfies {\it Condition $F^{++}$}. Below we denote $\sigma _i(A) :=
\{ \xi\in \R | \ i\xi \in \sigma (A)\}$.
\subsection{Proof of Theorem \ref{the 1}.} Since $u$ is a mild
solution of (\ref{eq}), by \cite[Lemma 4.7]{min} we have
\begin{equation}\label{est of spe of mil hl}
sp^+_{\cal F}u \subset \sigma_i(A)\subset \{ 0\} .
\end{equation}
Moreover, for $Re\lambda \not= 0$,
\begin{eqnarray}\label{4.12}
R(\lambda , \tD) \tilde{u} &=& R(\lambda ,\tilde{{\cal
A}})\tilde{u} .
\end{eqnarray}
Let $w={\cal A}u$. Then, by the assumption, $w\in BC(\R^+,\X)$, so
by the identity $R(\lambda ,\tD )\tD \tilde{u} = \lambda R(\lambda
,\tD )\tilde{u} -\tilde{u}$, we have
\begin{eqnarray}
R(\lambda ,\tD )\tilde{w} &=& R(\lambda ,\tD ){\tilde{\cal A}} \tilde{u} \nonumber\\
&=& R(\lambda ,\tD )\tD \tilde{u} \nonumber\\
&=&  \lambda R(\lambda ,\tD )\tilde{u} -\tilde{u}. \label{sp}
\end{eqnarray}
Therefore, by (\ref{est of spe of mil hl}),
\begin{equation}
sp^+_{\cal F}w \subset sp^+_{\cal F}u \subset \{ 0\} .
\end{equation}
Moreover, we have
\begin{eqnarray}
\lim_{\alpha \downarrow 0}   \| \alpha R(\alpha ,\tD )\tilde{w} \|
&=&
   \lim_{\alpha \downarrow 0} \| \alpha (\alpha  R(\alpha ,\tD )\tilde{u}
  -\tilde{u}) \|  \nonumber   \\
  &\le & \lim_{\alpha \downarrow 0} \| \alpha^2 R(\alpha ,\tD
  )\tilde{u}\|
+\lim_{\alpha \downarrow 0} \| \alpha \tilde{u}\|  \nonumber  \\
&\le & \lim_{\alpha \downarrow 0}\alpha^2 \| \frac{\tilde{u}}{\alpha
}\|  \nonumber   \\
&=& 0 .
\end{eqnarray}
Therefore, the assertion of Theorem \ref{the 1} follows from
Theorem \ref{the 3}.

\subsection{Proof of Theorem \ref{the 2}}
Let $\psi$ be any function in $L^\infty(\R^+, \X)$ (or
$L^\infty(\R^-, \X)$, respectively) recall that it can be
identified with its extension to $\psi \in L^\infty(\R, \X)$ by
setting $\psi (t)=0$ for all $t<0$ (or $t>0$, respectively).

\medskip
Theorem \ref{the 2} follows immediately from Theorem \ref{the 3}
and the following lemma
\begin{lemma}\label{lem 3}
Under the assumption of Theorem \ref{the 2},
\begin{equation}\label{empty}
sp^+_\cal F (w) =\emptyset ,
\end{equation}
where
\begin{equation}
w(t):= \int^\infty_0 \phi (s)u(t+s)ds , \mbox{for all t} \in \R^+.
\end{equation}
\end{lemma}
\begin{proof}
To prove Lemma \ref{lem 3}, we need several technical steps.
Let $\varphi (\theta ):= \phi (-\theta )$ for all $\theta \in \R$.
We re-write the function $w$ as
\begin{eqnarray*}
w(t) &=& \int^\infty_{-\infty} \phi (s)u(t+s)ds \\
&=&\int^\infty_{-\infty} \varphi (t-s)u(s)ds \\
&=& \varphi * u .
\end{eqnarray*}
First, we assume that the Fourier transform of $\phi$ (denoted by
$\F\phi$) vanishes in a neighborhood of $i\sigma (A) \cap \R $. As
shown below this assumption does not restrict the generality of
the proof. Since the Fourier transform $\F \varphi $ of $\varphi$
satisfies $\F \phi (- \xi )=\F \varphi (\xi )$ it vanishes in a
neighborhood of $-(i\sigma (A) \cap \R )= \R\backslash
\sigma_i(A)$, so, as is well known (see e.g. \cite[Chap. 4, and
5]{arebathieneu}), the Carleman transform $\widehat{\varphi}
(\lambda )$ of $\varphi$ is holomorphic in a neighborhood of
$\R\backslash \sigma_i(A)$. In turn, by the proof of
\cite[Proposition 2.3]{liunguminvu}, $ R(\lambda ,{\cal D})\varphi
$ should be holomorphic in a neighborhood of $\R\backslash
\sigma_i(A)$ as well. Next, by a simple computation, (by
considering $\varphi$ and $u$ as functions on $\R$) we can show
that
\begin{equation}\label{100}
R(\lambda ,{\cal D}) (\varphi *u) = R(\lambda ,{\cal D})\varphi *
u .
\end{equation}
 By Condition $F^{++}$, for each $g
\in \F$, $\psi \in L^1(\R^- )$, $\psi * g\in \F$. Therefore,
$\psi$ induces a map $BC(\R^+ ,\X)/\F\ni \tilde{g}\mapsto
\tilde{\psi} *\tilde{g} \in BC(\R^+ ,\X)/\F$. Notice also that if
$\psi_n \to \psi$ in $L^1(\R )$, then the induced maps
$\widetilde{\psi_n}\to \widetilde{\psi}$.   Next, consider the
canonical projection $p:BC(\R^+ ,\X)\to BC(\R^+ ,\X)/\F $ and the
restriction $r$ to the half line of elements of $BC(\R^+,\X)$.
Hence, from (\ref{100}),
\begin{equation}\label{3.9}
R(\lambda ,{\tD}) (\widetilde{\varphi *u}) =p\circ r(R(\lambda
,{\cal D}) (\varphi *u) = p \circ r(R(\lambda , {\cal D})\varphi
* u ).
\end{equation}
This shows that $R(\lambda ,{\tD}) (\widetilde{\varphi *u})) $ is
holomorphic at any $\xi\in \C$ whenever so is $ R(\lambda , {\cal
D})\varphi$. Now, by the above remark that $R(\lambda ,\cD
)\varphi$ is holomorphic in a neighborhood of $ \R\backslash
\sigma_i(A)$, (\ref{3.9}) yields in particular that
\begin{equation}\label{101}
sp^+_\F (\varphi *u) \subset \R \backslash \sigma_i (A) .
\end{equation}

\bigskip
On the other hand, a simple computation shows that
\begin{equation}\label{200}
R(\lambda ,{\cal D}) (\varphi *u) = \varphi * R(\lambda ,{\cal D})
u .
\end{equation}
Recall that we consider $u$ as a function on $\R$ with $u(t)=0$
for all $t<0$, so for $t\ge 0$ and $Re \lambda <0$,
$$
R(\lambda ,\cD )(\varphi *u )(t) = -\int^\infty _{-\infty} \varphi
(t-\xi )\int^\xi_0 e^{-\lambda (\zeta -\xi )}u(\zeta )d\zeta .
$$
Hence,
\begin{equation}\label{203}
R(\lambda ,\tD )(\widetilde{\varphi *u} ) = \tilde{\varphi }*
R(\lambda ,\tD )\tilde{u} .
\end{equation}
By \cite[Lemma 4.7]{min}, $sp^+_\F (u )\subset \sigma_i (A)$, $
R(\lambda ,\tD )\tilde{u}$ is holomorphic in a neighborhood of
$\sigma_i (A)$, so is $ R(\lambda ,\tD )(\widetilde{\varphi *u}
)$. In turn, this yields that
\begin{equation}\label{201}
sp^+_\F (\varphi *u) \subset \sigma_i (A).
\end{equation}
Now (\ref{101}) and (\ref{201}) prove (\ref{empty}).

\medskip
In the general case, if $\phi$ is approximated by $\phi_n$ (in
$L^1(\R)$) whose Fourier transform $\F\phi_n$ vanishes in a
neighborhood of $i\sigma (A)\cap \R$ we note that (\ref{100}) and
(\ref{203}) allows us to "pass to the limits". Therefore,
(\ref{empty}) is proved, so the proof of Theorem \ref{the 2} is
complete.
\end{proof}

\end{document}